\documentclass[letterpaper, 10 pt, conference]{ieeeconf}  

\IEEEoverridecommandlockouts                              

\usepackage{cite}
\usepackage{graphicx}
\usepackage{algorithm,algorithmic}
\usepackage{hyperref}
\usepackage{textcomp}

\let\labelindent\relax
\usepackage{enumitem}

\usepackage[utf8]{inputenc}
\usepackage[cmex10]{amsmath}
\usepackage{amsbsy,pifont,verbatim,amsfonts,amssymb}
\usepackage{array}
\usepackage[caption=false,font=normalsize,labelfont=sf,textfont=sf]{subfig}
\usepackage{textcomp}
\usepackage{stfloats}
\usepackage{url}
\usepackage{verbatim}
\usepackage{hyperref}
\usepackage{xcolor}
\usepackage{nomencl}
\usepackage{etoolbox}
\usepackage{float}
\newtheorem{theorem}{Theorem}
\newtheorem{lemma}[theorem]{Lemma}
\newcommand{\R}{\mathbb{R}}

\usepackage{algorithmic} 
\title{\LARGE \bf
Second-Order Optimization via Quiescence
}

\author{Aayushya Agarwal$^{1}$, Larry Pileggi$^{1}$ and,
Ronald Rohrer$^{1}$
\thanks{$^{1}$A.A., L.P., R.R. are with Electrical and Computer Engineering Department, Carnegie Mellon University, 5000 Forbes Ave, Pittsburgh, PA, 15213
        {\tt\small \{aayushya, pileggi, rr1j\}@andrew.cmu.edu}}%
}

\begin{document}

\maketitle
\thispagestyle{empty}
\pagestyle{empty}

\begin{abstract}

Second-order optimization methods exhibit fast convergence to critical points, however, in nonconvex optimization, these methods often require restrictive step-sizes to ensure a monotonically decreasing objective function. In the presence of highly nonlinear objective functions with large Lipschitz constants, increasingly small step-sizes become a bottleneck to fast convergence. We propose a second-order optimization method that utilizes a dynamic system model to represent the trajectory of optimization variables as an ODE. 
We then follow the quasi-steady state trajectory by forcing variables with the fastest rise time into a state known as \emph{quiescence}.
This optimization via quiescence allows us to adaptively select large step-sizes that sequentially follow each optimization variable to a quasi-steady state until all state variables reach the actual steady state, coinciding with the optimum. The result is a second-order method that utilizes large step-sizes and does not require a monotonically decreasing objective function to reach a critical point. Experimentally, we demonstrate the fast convergence of this approach for optimizing nonconvex problems in power systems and compare them to existing state-of-the-art second-order methods, including damped Newton-Raphson, BFGS, and SR1. 

\end{abstract}

\section{Introduction}
Second-order optimization methods have proven effective in solving large-scale problems. Methods such as Newton-Raphson and Quasi-Newton use curvature information to achieve faster convergence to critical points of the objective function, which is imperative for applications such as optimizing new infrastructure in power systems. 

Second-order methods, however, can often exhibit numerical oscillation or overflow without proper step-size damping \cite{royer2018complexity}. In practice, backtracking-style methods are frequently employed in conjunction with the optimization algorithm to select step-sizes that monotonically decrease the objective function at each iteration. Although this approach guarantees convergence, it often comes with the expense of small step-sizes that prove to be inefficient for highly nonconvex objective functions. One limitation of utilizing backtracking line search for high-dimensional problems stems from the multitude of directional axes. This complexity imposes constraints, as some axes necessitate exceedingly small step sizes, while others could easily accommodate larger ones, which are desirable for improving efficiency.  In practice, manual tuning of the step-sizes is often required to converge within a maximum number of allowable iterations.

In this work, we propose an adaptive, second-order optimization method that we call OptiQ, which was inspired by the circuit simulation algorithm (ACES) described in \cite{aces}. By leveraging a continuous time representation of the optimization variables, known as gradient flow \cite{behrman1998efficient}, we model the trajectory of optimization variables via state-space equations. The trajectory can then be simulated using numerical methods to obtain a steady state that coincides with the critical point of the objective function. OptiQ adaptively selects step-sizes that sequentially force each optimization variable to a quasi-steady state. Unlike traditional step-size selection methods that rely on monotonically decreasing the objective function, our proposed adaptive step-sizing can take large steps that are proportional to the remaining dominant time constant. The optimization variable with the smallest dominant time constant is then taken to quiescence, defined as the quasi-steady state of a single variable. Quiescent state variables are controlled by non-quiescent state variables (using a partial Hessian), and all non-quiescent variables follow the trajectory designated by the gradient flow. In this process, OptiQ sequentially takes each state variable to quiescence until all variables reach a steady state, which indicates arrival at the critical point of the objective function. The advantages of OptiQ are the following:
\begin{itemize}
    \item Adaptive step-size selection based on quasi-steady state behavior (quiescence) allows larger step sizes
    \item Speedup achieved by inverting only partial Hessians
\end{itemize}

In this work, we focus primarily on applications where second-order methods are preferable, such as power systems \cite{FOSTER2022108486}. Experimentally, we study the performance of this method in optimization test cases as well as in non-convex power system optimization problems to study infeasible planning and operation scenarios. Compared to other second-order adaptive optimization methods, OptiQ reaches a critical point in fewer iterations and less runtime.

\section{Related Work}

Most adaptive second-order optimization methods, such as Newton-Raphson, BFGS \cite{nocedal1999numerical}, and SR1 \cite{nocedal1999numerical}, use a backtracking line-search style approach to select a step-size that guarantees a monotonically decreasing objective \cite{truong2021new} (or Wolfe's conditions \cite{prudente2022quasi}). This approach has been used to solve large-scale optimization problems in power systems \cite{baker2020learning,giras1981quasi,tang2017real,semlyen2001quasi,domyshev2018optimal,alexander2020improved}, robotics \cite{kashyap2023dynamic,rehman2022solution,lee2005newton,jha2020quasi}, unconstrained optimization \cite{hassan2021adaptive,hassan2020new} and machine learning \cite{rafati2020quasi,yu2008quasi,byrd2016stochastic, bollapragada2018progressive, ma2020apollo, goldfarb2020practical}. However, this requirement creates highly restrictive step-sizes that are inefficient in the practical solution of nonconvex optimization problems.

In this work, we adopt a dynamical system perspective of the optimization problem, known as gradient flow \cite{behrman1998efficient}, \cite{attouch1996dynamical}, \cite{brown1989some}. The dynamical system model has been used to study and design first-order optimization methods \cite{helmke2012optimization},\cite{cortes2006finite}. Other works have used Lyapunov condition criteria to define the convergence rates of optimization methods \cite{wilson2018lyapunov}, \cite{wilson2021lyapunov}, \cite{polyak2017lyapunov}, \cite{hustig2019robust}, while others have designed predictor-corrector \cite{wilson2021lyapunov}, and proportional-integral control mechanisms to improve the convergence rates \cite{wadiaoptimization,maleki2021heunnet} of first-order optimization methods. 
Unlike prior works, OptiQ studies the quasi-steady state behavior of the gradient flow trajectory, circumventing the stability challenges of selecting an appropriate step-size in second-order optimization methods. 

\section{Problem Formulation}
In this work, we solve the following optimization problem:
\begin{subequations}
\label{eq:objective_fun}
    \begin{equation}
        \min_x f(x)
    \end{equation}
    \begin{equation}
        x^* \in \text{argmin}_x f(x)
    \end{equation}
\end{subequations}
where $x\in \R^n$ and the objective, $f:\R^n \rightarrow \R$, satisfies:
\begin{enumerate}[label=\textbf{(A\arabic*)}]
        \item $f\in C^2$ and $\inf_{x\in \R^n}f(x)>-R$ for some $R>0$. \label{a1}
        \item $f$ is second-order Lipschitz continuous \label{a4}
\end{enumerate}
The critical points of $f(x)$, defined as $x^*$ that ensure $\nabla f(x^* )=0$, lie within the set $S=\{x|\nabla f(x)=0\}$. 
\section{Dynamic System Model}
    The trajectory of the optimization variables is modeled using continuous-time gradient flow equations \cite{behrman1998efficient}:
    \begin{equation}
        \label{eq:gd_flow}
            \dot{x}(t) = -\nabla f(x(t)),\;\;
            x(0)= x_0.
    \end{equation}
    $x(t)$ models the trajectory of the optimization variable, $x^k$, and is initialized at $t=0$ in \eqref{eq:gd_flow} at an initial state of $x(t=0)=x_0\in\R^n$. For an objective function satisfying \ref{a1}-\ref{a4}, the gradient-flow reaches a steady state, $x^*$, when:
    \begin{equation}
        \dot{x}^*=-\nabla f(x^*)=0.
        \label{eq:steady_state_condition}
    \end{equation}
The state vector $x^* \in \R^n$ at the steady state is within the set of critical points of the objective function, $x^*\in S$.
\section{Explicit Numerical Integration}
\begin{figure}
    \centering
    \includegraphics[width=0.5\columnwidth]{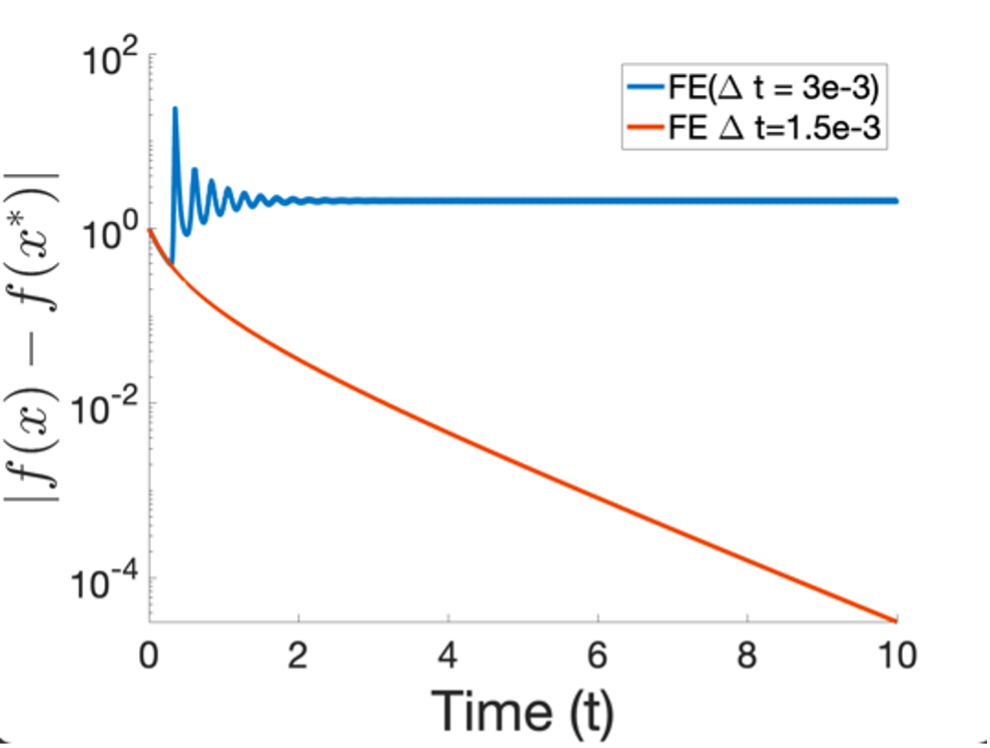}
    \caption{Solving a Rosenbrock function using an explicit numerical integration method (Matlab ODE23) results in numerical oscillations/divergence from the optimum. An FE with a smaller time-step (1.5e-3 s) converges to the local optimum but requires 5425 iterations (time-steps).}
    \label{fig:rosenbrock_ex}
\end{figure}

We achieve the steady-state solution of the gradient flow by marching through time in discrete steps until we reach convergence, defined by \eqref{eq:steady_state_condition}. The state at each time point, $x(t+\Delta t)$, is solved using the following: 
\begin{align}
    x(t+\Delta t) = x(t) + \int_t^{t+\Delta t} \dot{x}(t) dt \\
    x(t+\Delta t) = x(t) + \int_t^{t+\Delta t} -\nabla f(x(t)) dt.
    \label{eq:numerical_integration}
\end{align}
Generally, the integral on the right-hand side of \eqref{eq:numerical_integration} does not have an analytical solution and is instead approximated with the use of numerical methods. Numerical methods can be broadly classified as explicit or implicit integration. Explicit integration methods approximate the integral on the right-hand side of \eqref{eq:numerical_integration} using state information from previous time steps. The simplest (first order) explicit numerical integration method is Forward Euler (FE), which approximates \eqref{eq:numerical_integration} as:
\begin{align}
    x(t+\Delta t) &= x(t) + \Delta t \dot{x}(t) \\
    &= x(t) - \Delta t \nabla f(x(t)). \label{eq:fe_step}
\end{align}
FE is akin to gradient descent with a step-size of $\Delta t$ \cite{agarwal2023equivalent}. Explicit methods, such as Forward Euler, are known to suffer from numerical instability, where large time-steps can potentially diverge from the ODE trajectory or cause numerical instability that oscillates away from the critical point. In optimization, this is analogous to selecting large step-sizes that cause oscillations and divergence. An example of this is shown in Figure \ref{fig:rosenbrock_ex}, where a FE integration (i.e., gradient descent with a step size of $\Delta t$) is used to simulate the gradient flow response of optimizing a Rosenbrock function. A large time step causes divergence from the optimum.

From a numerical integration perspective, to guarantee stability, FE step-sizes must be restricted to the following:
\begin{equation}
    \Delta t \leq \frac{2}{\lambda},
    \label{eq:fe_bound}
\end{equation}
where $\lambda$ is the largest eigenvalue of the Hessian evaluated at time $t$ \cite{pillage1995electronic}. This time-step is especially restrictive in stiff gradient flow systems, making the method inefficient.

\section{Adaptive Numerical Integration via Quiescence}
Inspired by the circuit simulation algorithm ACES \cite{aces,esim}, we propose a new adaptive integration method that overcomes the  FE challenges of instability.
Rather than following the trajectory $\dot{x}(t)$, our approach designs a new trajectory $\dot{\Bar{x}}(t)$, which follows the quasi-steady state behavior. This trajectory assumes that small-time constants with low residue values of the dynamical system \eqref{eq:gd_flow} sequentially and instantaneously reach a quasi-steady state. This allows us to take larger steps that implicitly filter the small-time constants of the system \eqref{eq:gd_flow}, thereby avoiding the FE restriction in \eqref{eq:fe_bound}.
With our new algorithm, OptiQ, the new updates are now: 
\begin{equation}
    x(t+\Delta t) = x(t) + \Delta t \dot{\Bar{x}}(t).
    \label{eq:modified_ode}
\end{equation}
\subsection{Quiescent State Variables}
OptiQ's trajectory, $\dot{\Bar{x}}(t)$, sequentially follows each variable with the smallest dominant time-constant into a local quasi-steady state which we call \emph{quiescence}. 
A state-variable, $x_q$, is in quiescent when:
\begin{equation}
    \ddot{x}_q(t)=0. \label{eq:quiescence_def}
\end{equation}
 
Note that this definition of quiescence differs from that of steady state. While steady state is achieved when the time derivative of the entire state vector, $\dot{x}(t)$, is zero, quiescence is achieved at different times for each optimization variable. This approach does not follow the trajectory of gradient flow, but instead jumps to successive quiescent states. As the states with smaller dominant time-constants enter quiescence, OptiQ takes larger time-steps proportional to the time-constants of the slower state variables. 

At each iteration, the state vector is now composed of states in quiescence, $x_q\in \R^q$, and states not in quiescence, $x_{nq}\in \R^{n-q}$. By definition, non-quiescent states follow the original gradient flow trajectory, $\dot{x}_{nq} = -\frac{\partial f}{\partial x_{nq}}(t)$. The state vector, $x$, is then partitioned as:
\begin{equation}
    x = \begin{bmatrix}
        x_q \\ x_{nq}.
    \end{bmatrix}
    \label{eq:state_variable_vec}
\end{equation}
Because the second-order time derivative of all quiescent variables is zero, by definition, the trajectory of quiescent state variables is determined solely by non-quiescent variables. We study the second-order time derivatives of the state variables to design the updated trajectory, $\dot{\Bar{x}}(t)$. The second-order time derivative of the quiescent state-variables is:
\begin{subequations}
\begin{align}
    \ddot{x}_q(t) &= - \frac{d}{dt}\frac{\partial}{\partial x_q} f(x(t)) \\ 
    &= - \frac{\partial^2}{\partial x_q^2}f(x(t))\dot{\bar{x}}_q - \frac{\partial^2}{\partial x_q \partial x_{nq}}f(x(t)) \dot{x}_{nq}(t).
\end{align}    
\end{subequations}
Since the second-order time derivatives of the quiescent state variables are zero, we can solve for the time derivative of the quiescent state variables:
\begin{equation}
    \dot{\bar{x}}_q(t) = -\left( \frac{\partial^2}{\partial x_q^2}f(x(t)) \right)^{-1}\frac{\partial^2}{\partial x_q \partial x_{nq}}f(x(t)) \dot{x}_{nq}(t)
        \label{eq:quiescent_trajectory}
\end{equation}

The trajectory of the quiescent state variables is determined by the non-quiescent variables as well as the second-order terms $\frac{\partial^2}{\partial x_q^2}f(x(t))$ and $\frac{\partial^2}{\partial x_q \partial x_{nq}}f(x(t))$. Note, this update only requires inverting a subset of the Hessian,$\frac{\partial^2}{\partial x_q^2}f(x(t)) \in \R^{q\times q}$ , which reduces the complexity in comparison to second-order methods that require inverting the full Hessian (e.g., Newton-Raphson).

\subsection{Taking Variables to Quiescence}
	OptiQ progresses through time by forcing a set of non-quiescent variables with the fastest rise-time to quiescence at each iteration. The non-quiescent variable with the fastest rise time is determined by:
 \begin{equation}
     x_{nq_i} = \text{argmin}  \; \Tilde{\tau}_i \;\; \forall i \in \{n-q\},
     \label{eq:nq_to_q}
 \end{equation}
 where $\Tilde{\tau}_i$ is the time-constant associated to the first-order approximation of the non-quiescent state response and is calculated as:
 \begin{equation}
     \Tilde{\tau}_i = -\frac{\dot{x}_{nq_i}}{\ddot{x}_{nq_i}}
     \label{eq:time_constant_def}
 \end{equation}
Then, the entire system moves forward by a time-step of
\begin{equation}
    \Delta t = \min \Tilde{\tau}_i.
    \label{eq:delta_t_def}
\end{equation}
 
 To demonstrate that \eqref{eq:nq_to_q} is an appropriate choice for the non-quiescent variable that enters quiescence, we first prove that $-\frac{\dot{x}_{nq_i}(0)}{\ddot{x}_{nq_i}(0)}$ approximates the time-constant associated with the first-order approximation of the state-response. After establishing $\Tilde{\tau}_i = -\frac{\dot{x}_{nq_i}(0)}{\ddot{x}_{nq_i}(0)}$, we prove that $\text{argmin}  \Tilde{\tau}_i$ is the variable with the fastest rise-time that enters quiescence first.
 \begin{theorem}
     $-\frac{\dot{x}_{nq_i}(0)}{\ddot{x}_{nq_i}(0)}$ approximates the first-order time constant of each non-quiescent state variable, $x_{nq_i}$.
 \end{theorem}
 \begin{proof}
     Suppose the gradient flow equations are linearized around a point at time $T_0$ to produce the following linear set of ordinary differential equations (ODEs):
     \begin{equation}
         \dot{x}_{nq}(t) = Ax_{nq}(t) + Bu,
         \label{eq:linear_gd_flow}
     \end{equation}
     where $A$ represents the linearization of the Hessian, $\nabla^2 f(x(T_0))$, and $Bu$ is a step-input representing the gradient at time $T_0$ and is only active for $t\geq0$ and is zero for $t<0$.

     We define a shifted variable, $\hat{x}_{nq}(t) = x_{nq}(t) - x_{nq}(T_0)$, which is substituted into the linearized ODE as:
     \begin{align}
         \dot{\hat{x}}_{nq}(t) = A\hat{x}_{nq}(t) + \hat{B},
         \label{eq:shifted_linear_gd_flow}
     \end{align}
     where $\hat{B} = Bu + Ax_{nq}(T_0)$. The shifted version of the ODE \eqref{eq:shifted_linear_gd_flow} has an initial condition of $\hat{x}_{nq}(T_0)=0$.

    The response of the non-quiescent variables in the linearized system is defined by a sum of $n$ exponentials
    \begin{equation}
        \hat{x}_{nq_i}(t) = \gamma + \sum_{i=1}^n k_i e^{-p_i(t-T_0)},
    \end{equation}
    where $k_i$ and $p_i$ are the residue and pole of the linearized system \eqref{eq:shifted_linear_gd_flow} and $\gamma$ is a constant offset due to the input.

     We approximate the time-domain response by a first-order approximation, as often performed in circuit analysis \cite{awe}. The first-order approximation to the trajectory of each  variable, $\hat{x}_{nq_i}(t)$, is given as a single exponential response with a fitted pole, $\Bar{p}_i$ and residue, $\bar{k}_i$, as follows:
     \begin{equation}
         \hat{x}_{nq_i}(t) = \bar{k}_i e^{-\bar{p}_i (t-T_0)} + \gamma.
         \label{eq:approx_xt}
     \end{equation}
     The time derivative of the approximate trajectory is:
     \begin{equation}
         \dot{\hat{x}}_{nq_i}(t) = -\bar{k}_i \bar{p}_i e^{-\bar{p}_i (t-T_0)},
         \label{eq:approx_dxdt}
     \end{equation}
     and the second-order time derivative is
          \begin{equation}
         \ddot{\hat{x}}_{nq_i}(t) = \bar{k}_i \bar{p}^2_i e^{-\bar{p}_i (t-T_0)},
        \label{eq:approx_d2xdt2}
     \end{equation}
Dividing the two time-derivatives equals
\begin{equation}
    \frac{\dot{\hat{x}}_{nq_i}(t=T_0)}{\ddot{\hat{x}}_{nq_i}(t=T_0)} = -\frac{1}{\bar{p}_i} = -\Tilde{\tau}_i,
\end{equation}
     where $\Tilde{\tau}_i$ is the time-constant associated with the first-order approximation of the state trajectory.
 \end{proof}

 \begin{theorem}
     The non-quiescent variable with the smallest value of  $\Tilde{\tau}_i$ has the fastest rise time and is an appropriate choice to force to quiescence. 
 \end{theorem}
 \begin{proof}
     The response of the linearized ODE \eqref{eq:shifted_linear_gd_flow} can be analyzed in the Laplace domain as:
     \begin{equation}
         \hat{X}_{nq}(s) = A^{-1}(1-sA^{-1})^{-1}\frac{\hat{B}}{s},
     \end{equation}
     where $\hat{X}_{nq}(s)$ is the Laplace transformation of $\hat{x}_{nq}(t)$.
     The transfer function, $H(s)=\frac{X_{nq}(s)}{\hat{B}/s}$, is defined as
     \begin{equation}
         H(s) = A^{-1}(1-sA^{-1})^{-1},
                  \label{eq:transfer_function}
     \end{equation}
     and can be expanded as 
     \begin{equation}
         H(s)=\frac{1+a_1s+a_2 s^2 + \cdots + a_n s^n}{1+b_1s+b_2 s^2 + \cdots + b_m s^m},
         \label{eq:transfer_function_fraction}
     \end{equation}
     where $m>n$. Additionally, $H(s)$ can be defined as a multiplicative series of poles, $p_i$, and zeros, $z_i$, as follows:
     \begin{equation}
         H(s) = K\frac{(1-\frac{s}{z_1})(1-\frac{s}{z_2})\cdots(1-\frac{s}{z_n})}{(1-\frac{s}{p_1})(1-\frac{s}{p_2})\cdots(1-\frac{s}{p_m})}.
                  \label{eq:transfer_function_multiples}
     \end{equation}
     From analyzing \eqref{eq:transfer_function_fraction} and \eqref{eq:transfer_function_multiples}, we conclude that:
     \begin{equation}
         a_1 = \sum_{j=1}^n -\frac{1}{z_j} \;\;\;, b_1 = \sum_{j=1}^m -\frac{1}{p_j}.
         \label{eq:definition_a1b1}
     \end{equation}
     Expanding about $s=0$ \cite{pillage1995electronic}, the transfer function can be defined as a sum of moments:
     \begin{equation}
         H(s) = m_0 + m_1 s + m_2 s^2 + \cdots ,
     \end{equation}
     where $m_0=K$ and $m_1 = a_1 - b_1$. 
 In the first-order approximation of the trajectory \eqref{eq:approx_xt}, the transfer function is approximated by the first two moments \cite{awe}:
     \begin{equation}
         H(s) \approx m_0 + m_1 s = \frac{\bar{k}}{s-\bar{p}},
     \end{equation}
     where $\bar{k},\bar{p}$, are the identical residue/pole pair as in \eqref{eq:approx_xt}.

     Fitting the two moments to the poles and zeros of the original transfer function, \eqref{eq:transfer_function}, we determine that
     \begin{align}
         m_0 = \frac{\bar{k}}{\bar{p}}, \;\;
         m_1 = \frac{\bar{k}}{\bar{p}^2}.
     \end{align}
     This means that
     \begin{align}
         \frac{m_0}{m_1} = \bar{p} &= \frac{1}{a_1 - b_1} =  \frac{1}{\sum_{j=1}^m \frac{1}{p_j} - \sum_{j=1}^n \frac{1}{z_j}},
     \end{align}
     and the first-order approximate time constant, $\Tilde{\tau}_i=1/\bar{p}_i$, is
     \begin{equation}
         \Tilde{\tau}_i=\sum_{j=1}^m \frac{1}{p_j} - \sum_{j=1}^n \frac{1}{z_j}.
         \label{eq:time_constant_poles_zeros}
     \end{equation}
     Note, all variables share the same poles, $p_i$, since the poles are derived from the eigenvalues of the system matrix, $A$. This implies the variable with the smallest approximate time-constant in \eqref{eq:time_constant_poles_zeros} will have the largest sum of zeros:
     \begin{equation}
         \text{argmin}_i \Tilde{\tau}_i= \text{argmax}_i \sum_{j=1}^n \frac{1}{z_j}.
     \end{equation}
     Although the zeros in the transfer function do not affect the settling time of a step response, they do decrease the rise time \cite{franklin2002feedback}. As a result, the variable with the maximum sum of zeros will have the fastest rise time and enter quiescence first.
    This is observed in Figure \ref{fig:quadratic_example}, where we optimize a quadratic objective with two variables. Both variables exhibit an identical settling time (due to a common set of poles), however, $x_1$ clearly has a faster rise time since $x_1$ has a larger sum of zeros compared to that of $x_2$. Therefore, $\text{argmin}_i \Tilde{\tau}_i$ provides an appropriate choice for $x_1$ to enter quiescence.
 \end{proof}

At each iteration, we force the non-quiescent variable with the smallest dominant time-constant into quiescence by taking a Forward-Euler time-step equal to $\Delta t = \min \Tilde{\tau}_i$. A step size proportional to the smallest time constant avoids the issue of small step sizes by implicitly assuming that variables associated with any smaller time constants have already reached a quasi-steady state. 
\begin{theorem}
    A time-step of $\Delta t = \min \Tilde{\tau}_i$ is bounded by the smallest and largest time-constants of the linearized system matrix, $A$ from \eqref{eq:shifted_linear_gd_flow}.
\end{theorem}
\begin{proof}
    See Appendices \ref{sec:appendix_lower_bound}, \ref{sec:appendix_upper_bound}.
\end{proof}

 \subsection{Moving States Out of Quiscence}
State variables may fall out of quiescence if the quiescent approximation of $\ddot{x}_q=0$ is violated during the optimization trajectory. This occurs when the linear assumption of the quiescent trajectory in \eqref{eq:quiescent_trajectory} no longer hold due to the nonlinearities of the true gradient flow \eqref{eq:gd_flow}. To account for states moving out of quiescence, we define an error function that checks whether the quiescent state variables are deviating from the intended gradient flow trajectory. A vector of error values, $err(x)$, is calculated as:
\begin{equation}
    err(x) = | \dot{\bar{x}}_q(t+\Delta t) + \frac{\partial}{\partial x_q}f(x(t+\Delta t)) |,
\end{equation}
where $\dot{x}_q$ is evaluated using a difference of terms:
\begin{equation}
    \dot{\bar{x}}_q(t+\Delta t) = \frac{1}{\Delta t} (x_q(t+\Delta t) - x_q(t)).
\end{equation}
Each element in the vector, $err$, corresponds to the error in the quiescent trajectory in \eqref{eq:quiescent_trajectory}. A zero error means that the quiescent state follows the designated trajectory of $\ddot{x}_q=0$. 

However, the explicit nature of FE makes it prone to numerical instability whereby errors accumulate, leading to divergence from the steady-state (i.e., the local optimum). To guarantee that after $N$ iterations we have not diverged from the steady state (numerically defined by $\|f(x)\|<\eta$, where $\eta>0$), we move a state out of quiescence when:
\begin{equation}
    \max(err(x)) > \eta/N,
\end{equation}
where $\eta$ is a predefined scalar representing the tolerance for convergence for the optimization problem ($\| \nabla f(x) \| < \eta$) and $N$ is the maximum number of allowable iterations. This condition bounds the error at each iteration to ensure the worst-case error after $N$ iterations is $\eta$.

\subsection{Quasi-Steady State Gradient Flow Algorithm}
	The OptiQ algorithm, shown in Algorithm \ref{optiq_algorithm}, follows variables into quiescence until it reaches a critical point of the objective. At each iteration, we take step sizes equal to the estimated dominant time constant, avoiding the requirement of a monotonically decreasing objective function. 
\begin{algorithm}
\caption{OptiQ Algorithm}
\label{optiq_algorithm}
\textbf{Input:} $f(\cdot), x_0, \eta, N$
\begin{algorithmic}[1]
\STATE{$x = x_0, x_{nq} = x_0, x_q = []$}
\STATE{\textbf{do while} $\|\nabla f(x)\|^2 > \eta $}
\STATE{\hspace*{\algorithmicindent} $\dot{x}_{nq}(t)=\frac{\partial}{\partial x_{nq}} f(x(t)), \ddot{x}_{nq}(t)=\frac{\partial^2}{\partial x_{nq}^2} f(x(t)) \dot{x}_{nq}(t)$}
\STATE{\hspace*{\algorithmicindent} $[\Tilde{\tau}] = [-\frac{\dot{x}_{nq_i}}{\ddot{x}_{nq_i}}] \forall i \in \{n-q\}$}
\STATE{\hspace*{\algorithmicindent} $\Delta t = min([\Tilde{\tau}])$}
\STATE{\hspace*{\algorithmicindent} Append state variable with $argmin([\Tilde{\tau}])$ to $x_q$}
\STATE{\hspace*{\algorithmicindent} $\dot{x}_q(t)= - (\frac{\partial^2}{\partial x_q^2} f(x))^{-1} \frac{\partial^2}{\partial x_q \partial x_{nq}}f(x) \dot{x}_{nq}(t)$}
\STATE{\hspace*{\algorithmicindent} $x_{nq}(t+\Delta t) = x_{nq}(t) + \Delta t \dot{x}_{nq}(t)$}
\STATE{\hspace*{\algorithmicindent} $x_{q}(t+\Delta t) = x_{q}(t) + \Delta t \dot{x}_{q}(t)$}
\STATE{\hspace*{\algorithmicindent} $err(x) = | \dot{x}_q(t+\Delta t) + \frac{\partial}{\partial x_q}f(x(t+\Delta t)) |$}
\STATE{\hspace*{\algorithmicindent} \textbf{if } $err(x_i)>\eta/N$: set $x_{q_i}$ to non-quiescent state}
\end{algorithmic}
\end{algorithm}

\section{Convergence for Convex Functions}
We analyze the convergence of OptiQ’s quasi-steady state trajectory for convex functions using a positive definite Lyapunov function, $V(x)$. 
\begin{theorem}
    At each iteration of OptiQ, the quiescent trajectory asymptotically converges towards the critical point.
\end{theorem}
\begin{proof} To analyze the convergence of OptiQ, we define the following Lyapunov function to study \eqref{eq:modified_ode},
\begin{equation}
    V(\dot{x}) = \frac{1}{2} \| \dot{x}(t) \|^2.
    \label{eq:lyapunov_function}
\end{equation}
As observed by the second norm, $V(\dot{x})>0 \forall x \in \R^n \notin \{0\}$ and $V(0)=0$, \eqref{eq:lyapunov_function} is a positive function. By the chain rule, the time-derivative of the Lyapunov function is defined as:
\begin{align}
    \frac{d}{dt}V(\dot{x}) &= \frac{d}{dt}\frac{1}{2} \| \dot{x} \|^2 \\
    &= \dot{x}^T\ddot{x}(t).
    \label{eq:lyapunov_derivative}
\end{align}
With the state-vector defined as \eqref{eq:state_variable_vec}, eq \eqref{eq:lyapunov_derivative} is defined as:
\begin{equation}
    \frac{d}{dt} V(\dot{x}) = \begin{bmatrix} \dot{x}_q \\ \dot{x}_{nq}\end{bmatrix} ^T\begin{bmatrix} \ddot{x}_q \\ \ddot{x}_{nq}\end{bmatrix}.
\end{equation}
By defintion of quiescence in \eqref{eq:quiescence_def}, $\ddot{x}_q=0$, and the time-derivative of the Lyapunov function reduces to:
\begin{equation}
    \frac{d}{dt} V(\dot{x}) = \dot{x}_{nq}^T \ddot{x}_{nq}.
\end{equation}
This implies that non-quiescent variables solely dictate the convergence to a critical point. The second-order time derivative of the non-quiescent variables is defined as:
\begin{align}
    \ddot{x}_{nq} &= -\frac{d}{dt}\frac{\partial}{\partial x_{nq}} f(x) = - \frac{\partial^2}{\partial x_{nq}^2}f(x) \dot{x}_{nq}.
\end{align}
Therefore, the time-derivative of the Lyapunov function is
\begin{equation}
    \frac{d}{dt} V(\dot{x}) = -\dot{x}_{nq}^T\frac{\partial^2}{\partial x_{nq}^2}f(x),
\end{equation}
where $\frac{\partial^2}{\partial x_{nq}^2}f(x)$ is a principle sub-matrix of the Hessian, $\nabla^2 f(x)$. The convergence proof of a strictly convex function, $f(x)$, relies on the following lemma.
\begin{lemma}
    \label{lemma_convex}
    For a strictly convex function, $f(x)$, the Hessian is positive-definite, $\nabla^2 f(x) \succ 0$, and the principle sub-matrix is also positive-definite, $\frac{\partial^2}{\partial x_{nq}^2}f(x) \succ 0$ .
\end{lemma}
\begin{proof}
    See Observation 7.1.2 from \cite{horn2012matrix}.
\end{proof}
As a consequence of Lemma \ref{lemma_convex}, the principle sub-matrix, $\frac{\partial^2}{\partial x_{nq}^2}f(x)$ is positive definite, and can conclude the following
\begin{equation}
    \frac{d}{dt} V(\dot{x}) = -\dot{x}_{nq}^T\frac{\partial^2}{\partial x_{nq}^2}f(x) < 0 \forall \{\dot{x}\in \R^n|\dot{x}\neq 0\},
\end{equation}
which satisfies the last condition for Lyapunov stability and demonstrates that the quasi-steady state trajectory designed by OptiQ is asymptotically convergent to a critical point. Additionally, the convergence of OptiQ in \eqref{eq:quiescent_trajectory} is only dictated by the non-quiescent variables, $x_{nq}$, as the time-derivative of the proposed Lyapunov function, $d/dt V(\dot{x})$, reduces to a function of non-quiescent trajectories. This implies that as variables enter quiescence at each iteration, the convergence rate is determined by the slower active variables, $x_{nq}$.
\end{proof}
\subsection{Quadratic Objective Example}
We demonstrate the benefit of OptiQ for the optimization of the following quadratic objective function:
\begin{equation}
    \min_{x_1,x_2} 0.5(x_1 - 1)^2 + 50(x_1 - x_2)^2.
\end{equation}
The gradient-flow equations are:
\begin{align}
    \dot{x}_1 = -(x_1 - 1) - 100 (x_1 - x_2)\\
    \dot{x}_2 = 100(x_1-x_2).
\end{align}

\begin{figure}
        \centering
        \includegraphics[width=0.8\linewidth]{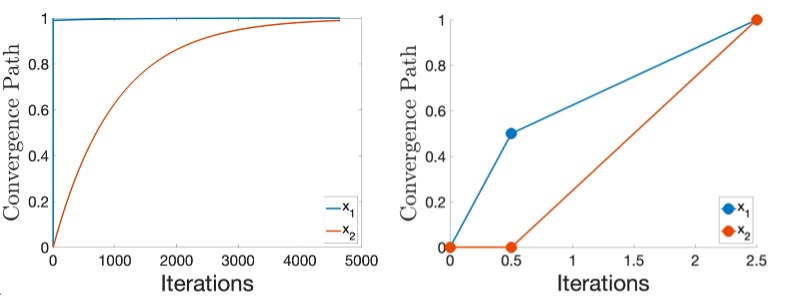}
        \caption{Transient response of gradient-flow using a) FE with time-step of 1s and b) OptiQ}
        \label{fig:quadratic_example}
\end{figure}

Solving this gradient-flow problem using FE integration with step-sizes that adhere to the specifications in \eqref{eq:fe_bound} requires over 5000 iterations to reach the optimum solution. The trajectories of $x_1$ and $x_2$, obtained through FE integration with a fixed time-step of 1, are depicted in Figure \ref{fig:quadratic_example}. As observed in Figure \ref{fig:quadratic_example}, $x_1$ reaches its steady-state much faster than $x_2$. This depicts the benefit of the quiescence concept from \cite{aces}. OptiQ can exploit the quasi-steady state behavior to force $x_1$ to a state of quiescence in the first iteration, while $x_2$ remains dormant. In the next iteration, OptiQ takes a step to force $x_2$  to quiescence, with the trajectory of $x_1$ determined by the dynamics of $x_2$. This process effectively avoids the issue of the small time-constant associated with $x_1$ and reaches the optimum in 2 iterations. 

\section{Results}
We evaluate the performance of OptiQ by solving convex and nonconvex optimization functions. We demonstrate the following benefits of our approach:
\begin{enumerate}
    \item Large, adaptive step sizes that do not require inner loops to satisfy a convergence condition
    \item Inverting a partial Hessian reduces the overall wall-clock time in comparison to NR
\end{enumerate}
We apply OptiQ to optimization test functions as well as solving large power systems optimizations, benchmarking it against state-of-the-art adaptive second-order methods, including damped Newton-Raphson, BFGS, and SR1. Step sizes for the three comparison methods are determined using a back-tracking line search with the Armijo condition \cite{armijo1966minimization}. 
\subsection{Optimization Test Functions}
\begin{figure*}
    \centering
\includegraphics[width=0.9\textwidth]{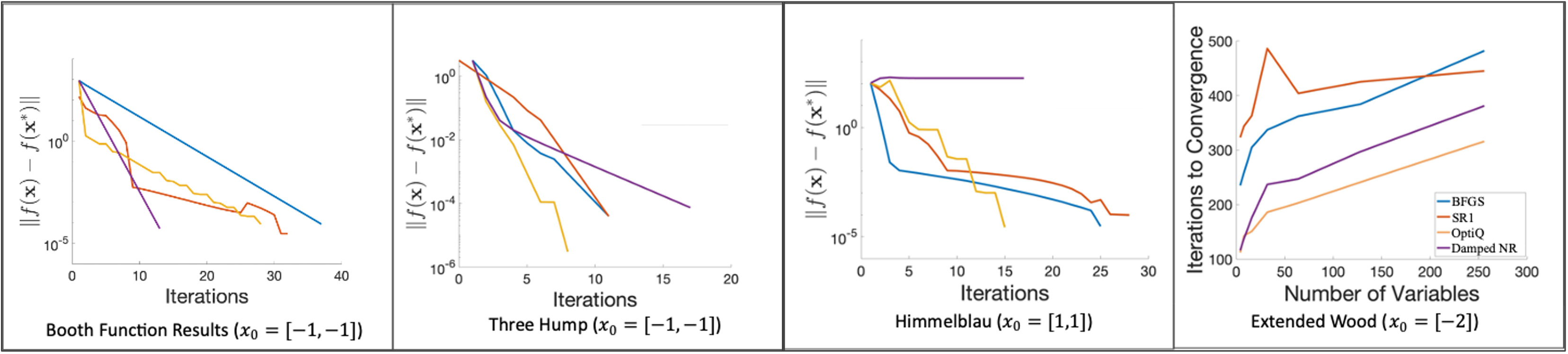}
    \caption{Comparison of second-order methods in optimizing test functions}
    \label{fig:test_functions}
\end{figure*}

The following convex and nonconvex test functions provide varying degrees of nonconvexities to test our optimization approach. The selected testcases are the convex Booth function \cite{testfunc}, the nonconvex Three Hump function, the nonconvex Himmelblau function \cite{himmelblau}, and the nonconvex Extended-Wood function \cite{testfunc}. To demonstrate scalability, we extend the Extended-Wood function to $n=256$ variables. The results for these functions are shown in Figure \ref{fig:test_functions}.

We observe that damped NR converges faster than OptiQ for the convex Booth function, as the convex solution space enables larger step sizes and reduces the need for multiple backtracking line search iterations.
 However, for optimizing nonconvex testcases, the comparison methods required multiple iterations of backtracking line search to find a satisfactory step size. On the other hand, OptiQ adaptively \emph{increases} the time steps (step sizes) during the solution process (a common restriction of backtracking line search methods \cite{armijo1966minimization}) to follow the quasi-steady state trajectory. This improves the rate of convergence, as it does not restrict the maximum value of the step size to a value of one, as is common in optimization methods.

The benefits of following the quasi steady-state behavior are greatly demonstrated in scaling the Extended wood function, which encompasses multiple repeated poles as the function is expanded to $n=256$ variables. OptiQ leverages such functions by making multiple state variables achieve quiescence simultaneously, resulting in fewer iterations required, as shown in Figure \ref{fig:test_functions}.

\subsection{Optimizing Power Systems}
The efficiency of OptiQ is next demonstrated for a power grid optimization problem, where we analyze the feasibility of a power grid configuration by optimizing the following:
\begin{equation}
    \min_x \| f(x) \|^2.
    \label{eq:power_systems_opt}
\end{equation}
$f(x):\R^n \rightarrow \R^n$ represents the network constraints of the grid, and $x\in \R^n$ is the state vector of bus voltages and reactive powers. This optimization problem tests the feasibility of networks and locates areas where new devices can be potentially added to improve stability. The objective \eqref{eq:power_systems_opt}, the details of which are provided in \cite{FOSTER2022108486}, is a multimodal, nonconvex function, that is challenging to solve for large, stressed systems. In this experiment, we study the feasibility of three stressed networks (IEEE-14, IEEE-500, and Pegase-13529 bus) under increased load conditions using OptiQ, Newton-Raphson, BFGS, and SR1. All optimizers reach the same optimum, with iteration counts shown in Table \ref{tab:power_systems_iterations}.

\begin{table}
    \centering
     \caption{Iterations performed for optimizing power system networks for infeasibility using BFGS, SR1, NR and OptiQ}
   \begin{tabular}{|c|c|c|c|}
    \hline
        Iterations & 14-Bus & 500-Bus & 13,529-Bus\\
        \hline
        BFGS & 7 & 14 & 116\\
        \hline
         SR1& 14 & 22 &68 \\
         \hline
         NR& 5 & 6 & 7\\
         \hline
         OptiQ & 4 & 4 &5 \\
         \hline
    \end{tabular}
    \label{tab:power_systems_iterations}
\end{table}

\begin{table}
    \centering
     \caption{Runtime normalized to NR for optimizing power system networks for infeasibility using BFGS, SR1, NR and OptiQ}
   \begin{tabular}{|c|c|c|c|}
    \hline
        Normalized Runtime & 14-Bus & 500-Bus & 13,529-Bus\\
            \hline
        BFGS & 1.18 & 1.86 & 10.61\\
             \hline
SR1& 2.14 & 3.08 & 6.05 \\
             \hline
NR& 1 & 1 & 1\\
             \hline
OptiQ & 0.91 & 0.74 & 0.68 \\
    \hline
    \end{tabular}
    \label{tab:power_systems_runtime}
\end{table}

	As shown in Tables \ref{tab:power_systems_iterations} and \ref{tab:power_systems_runtime}, OptiQ finds a solution in fewer iterations and less runtime than NR, BFGS, or SR1. Unlike adaptive line search methods \cite{armijo1966minimization}, OptiQ does not require multiple iterations to find a valid step-size. This reduces the per-iteration and overall wall-clock time. Also, unlike Newton-Raphson, OptiQ avoids inverting the full Hessian at each iteration. Since inverting a sparse Hessian is generally $O(n^{\sim 1.6} )$ \cite{pillage1995electronic}, reducing the size of the factored portion of the Hessian results in immediate savings.

\section{Conclusion}
We introduced a new second-order optimization method, OptiQ, which uses a dynamical system model of optimization to follow the quasi-steady state trajectory of optimization variables towards a critical point of the objective function. This process allows us to take large step-sizes that approximate the dominant time constant, thereby obviating the need for selecting a step size that monotonically decreases the objective function. This results in a fast, convergent second-order optimization method that: 1) adaptively selects large step-sizes; and 2) only requires a partial Hessian inversion at each iteration. We demonstrated that this method reduces the number of iterations as well as wall-clock time in optimizing both small nonconvex optimization test cases and large practical power systems examples. 
\appendix
\section{Appendix}

\subsection{Lower-Bound of Step-size}
\label{sec:appendix_lower_bound}
The lower-bound of the step-size, $\Delta t = \min \Tilde{\tau}_i$ in \eqref{eq:delta_t_def}, is $\Delta t \geq 1/p_1$, where $p_1$ is the largest pole in the linearized system matrix, $A$, in \eqref{eq:shifted_linear_gd_flow}.

\begin{proof}
 The trajectory of the non-quiescent state variables in the linearized ODE \eqref{eq:shifted_linear_gd_flow} is given by the following sum of exponentials:
 \begin{equation}
     x_{nq}(t) = \gamma + \sum_{i=1}^n k_i e^{-p_i (t-T_0)},
 \end{equation}
 where $k_i,p_i>0$ and real for a positive semi-definite, matrix, $A$, as demonstrated by Foster’s reactance theorem \cite{van1960introduction}. 

The first and second time-derivatives for the response are:
\begin{align}
    \dot{x}_{nq}(t) = \sum_{i=1}^{n} -k_i p_i e^{-p_i (t-T_0)}\\
    \ddot{x}_{nq}(t) = \sum_{i=1}^{n} k_i p_i^2 e^{-p_i (t-T_0)}.
\end{align}

Evaluated at $t=T_0$, we observe that the element-wise division of the first and second time-derivatives is as follows:
\begin{equation}
    -\frac{\dot{x}_{nq_i}(T_0)}{\ddot{x}_{nq_i}(T_0)} = \frac{\sum k_i p_i}{\sum k_i p_i^2}.
\end{equation}
Let $p_1$ be the largest pole for the system (with an associated residue $k_1$), then $\frac{\dot{x}_{nq_i}(T_0)}{\ddot{x}_{nq_i}(T_0)}$ can be expressed as
\begin{align}
    -\frac{\dot{x}_{nq_i}(T_0)}{\ddot{x}_{nq_i}(T_0)} &= \frac{p_1}{p_1^2}\left ( \frac{k_1 + \sum_{i>1}k_ip_i/p_1}{\sum_{i>1}k_ip_i^2/p_1^2}    \right ) \\
    &= \frac{1}{p_1}\left ( \frac{k_1 + \sum_{i>1}k_ip_i/p_1}{\sum_{i>1}k_ip_i^2/p_1^2}  \right )
\end{align}

Since $p_1 > p_i \;\; \forall i>1$, then
\begin{equation}
    \frac{k_1 + \sum_{i>1}k_ip_i/p_1}{\sum_{i>1}k_ip_i^2/p_1^2} \geq 1,
\end{equation}
which implies that
\begin{equation}
    -\frac{\dot{x}_{nq_i}(T_0)}{\ddot{x}_{nq_i}(T_0)} = \Tilde{\tau}_i \geq \frac{1}{p_1}.
\end{equation}
This proves that the time-step, $\Delta t=min \Tilde{\tau}_i$, is lower-bounded by the smallest time-constant of the system, $1/p_1$.

 \end{proof}

 \subsection{Upper-bound of Step-Size}
 \label{sec:appendix_upper_bound}

To determine the upper-bound of $\Delta t = \min \Tilde{\tau}_i$, we study the linearized system in \eqref{eq:shifted_linear_gd_flow} which exhibits multiple poles that are close to the smallest pole, $p_d$, defined as:
\begin{equation}
    p_d < p_i \forall i \neq d.
\end{equation}
Suppose a subset of the poles, $p_j\in p$ are within a small $\epsilon$-ball away from $p_d$, defined as
\begin{equation}
    p_j = p_d + \epsilon_j,
\end{equation}
where $0<\epsilon_j \leq \epsilon$. Then we can approximate the response of the non-quiescent variables by the subset of poles as:
\begin{equation}
    x_{nq}(t) \approx \gamma + k_de^{-p_d (t-T_0)} + \sum_j k_je^{-(p_d+\epsilon_j)(t-T_0)},
\end{equation}
where $k_d,k_j>0$ are the residues corresponding to the dominant pole, $p_d>0$ and poles, $p_j>0$, respectively.

The first and second time-derivatives of the non-quiescent state variables are:
\begin{equation}
    \dot{x}_{nq_i}(t) = -k_dp_de^{-p_d(t-T_0)} - \sum_j k_j (p_d + \epsilon_j)e^{-(p_d+\epsilon_j)(t-T_0)}
\end{equation}
\begin{equation}
    \ddot{x}_{nq_i}(t) = k_dp_d^2 e^{-p_d(t-T_0)} + \sum_j k_j (p_d + \epsilon_j)^2 e^{-(p_d+\epsilon_j)(t-T_0)}
\end{equation}
The ratio of the time-derivatives at $t=T_0$ is then defined as
\begin{align}
    \frac{\dot{x}_{nq_i}(T_0)}{\ddot{x}_{nq_i}(T_0)} &= - \frac{k_dp_d + \sum_j k_j (p_d + \epsilon_j)}{k_dp_d^2 + \sum_j k_j (p_d + \epsilon_j)^2} \\
    &= - \frac{k_dp_d + \sum_j k_j (p_d + \epsilon_j)}{k_dp_d^2 + \sum_j k_j (p_d^2 + 2p_d\epsilon_j + \epsilon_j^2)}
\end{align}
Assuming a small value of $\epsilon$, we approximate the ratio as
\begin{align}
    \frac{\dot{x}_{nq_i}(T_0)}{\ddot{x}_{nq_i}(T_0)} &\approx - \frac{k_dp_d + \sum_j k_j (p_d + \epsilon_j)}{k_dp_d^2 + \sum_j k_j (p_d^2 + 2p_d\epsilon_j)}\\
    &= -\frac{p_d(k_d + \sum_j k_j) + \sum_j k_j \epsilon_j}{p_d^2(k_d+\sum_j k_j) + \sum_j k_j2p_d\epsilon_j} \\
    & = -\frac{1}{p_d}\frac{(k_d + \sum_j k_j) + \sum_j k_j \epsilon_j/p_d}{(k_d+\sum_j k_j) + \sum_j k_j2\epsilon_j/p_d}.
\end{align}
Since we have defined $\epsilon>0$, this implies that
\begin{equation}
    \frac{(k_d + \sum_j k_j) + \sum_j k_j \epsilon_j/p_d}{(k_d+\sum_j k_j) + \sum_j k_j2\epsilon_j/p_d} <1.
    \label{eq:appendix_ineq}
\end{equation}
Furthermore, we can now bound $\frac{\dot{x}_{nq_i}(T_0)}{\ddot{x}_{nq_i}(T_0)}$ as
\begin{equation}
    \frac{\dot{x}_{nq_i}(T_0)}{\ddot{x}_{nq_i}(T_0)} < -\frac{1}{p_d},
\end{equation}
and therefore
\begin{equation}
    \Delta t = \min \left( -\frac{\dot{x}_{nq_i}(T_0)}{\ddot{x}_{nq_i}(T_0)} \right) < \frac{1}{p_d}.
\end{equation}

Additionally, when we have a repeated dominant pole, then $\epsilon_j \rightarrow 0$, which implies that the inequality in \eqref{eq:appendix_ineq} is now
\begin{equation}
       \frac{(k_d + \sum_j k_j) + \sum_j k_j \epsilon_j/p_d}{(k_d+\sum_j k_j) + \sum_j k_j2\epsilon_j/p_d} \rightarrow 1,
\end{equation}
and therefore 
\begin{equation}
    \Delta t = \min \left( -\frac{\dot{x}_{nq_i}(T_0)}{\ddot{x}_{nq_i}(T_0)} \right)\rightarrow \frac{1}{p_d} .
\end{equation}
This shows that the step-size defined by \eqref{eq:delta_t_def} is upper-bounded by the largest time-constant of the system, $1/p_d$.
\bibliographystyle{IEEEtran}
\bibliography{IEEEabrv,myBib}

\end{document}